\documentclass[a4paper,11pt]{amsart}
\usepackage{amsmath}
\usepackage{amssymb}
\usepackage{amsfonts}

\def \endproof {\quad \hfill  \rule{2mm}{2mm} \par\medskip}

\DeclareMathSymbol{\subsetneqq}{\mathbin}{AMSb}{36}

\def \endproof {\quad \hfill  \rule{2mm}{2mm} \par\medskip}
\DeclareMathSymbol{\subsetneqq}{\mathbin}{AMSb}{36}

\newcommand{\R}{\mathbb{R}}

\newtheorem{th1}{{\bf Theorem}}[section]

\newtheorem{lem}{{\bf Lemma}}[section]

\theoremstyle{remark}
\newtheorem{rem}{\bf Remark}[section]

\theoremstyle{definition}

\author{\bf Jamel BENAMEUR}
\address{Facult\'e des Sciences de Bizerte, D\'epartement
de Math\'ematiques\\ Zarzouna 7021, Tunisia.} \email{\it
jamel.benameur@fsb.rnu.tn}

\title{SMOOTHING EFFECTS FOR  NAVIER-STOKES EQUATIONS }
\date{\today}
\begin{document}
\begin{abstract}
We prove some smoothing effects for the 3-D Navier-Stokes equations for initial data belonging to the critical Sobolev space $H^{1/2}(\R^3)$. Asymptotic behavior of the global solution when the time goes to infinity is studied. We also obtain a new energy estimate. Other results in this direction and with different methods can be found in \cite{C4}.
\end{abstract}
\maketitle

\section{Introduction}
\noindent The purpose of this text is to establish some regularity
results for the 3-D incompressible Navier-Stokes equations on the whole space $\mathbb R^3$. Throughout this paper we consider the three-dimensional
incompressible Navier-Stokes equations\\
$$ \left\{\begin{array}{l}
  \displaystyle\partial_t
u-  \nu\Delta u+ (u.\nabla) u=-\nabla p,\quad
\mbox{on}\;\; \R_+\times\mathbb R^3, \\\\
div\; (u) = 0 \quad
\mbox{on}\; \; \R_+\times\mathbb R^3,  \\\\
u _{\mid t=0}=u^0 \quad\mbox{on}\;\;\mathbb R^3 ,\\
\end{array}\right.
\leqno(NS_\nu) $$ where $\nu>0$ is the viscosity of the fluid, and $u=u(t,x)=(u_1,u_2,u_3)$
and  $p=p(t,x)$ denote respectively the unknown velocity and
the unknown pressure of the fluid at the point
$(t,x)\in\R_+\times\mathbb R^3$. While
$u^0=(u^0_1(x),u^0_2(x),u^0_3(x))$ is a given initial velocity. If $u^0$ is quite regular, the divergence free condition
determine the pressure $p$. Moreover, $p$ can be expressed as
follows
$$p=-\Delta^{-1}\sum_{j,k}\partial_j\partial_k(u_ju_k).$$
The above problem has been studied by many authors like  \cite{CG},
\cite{FK},\cite{Le1}...  Using compactness methods, Leray proved In
1934 (\cite{Le1})   for $u^0\in L^{2}(\R^3)$ an existence result in
$L^\infty(\R_+,L^2(\R^3))\cap L^2(\R_+,\dot H^1(\R^3))$ for the
problem $(NS_\nu)$. Also, in two dimension space, Leray proved the
existence and uniqueness in ${\mathcal C}_b(\R_+,L^2(\R^2))\cap
L^2(\R_+,\dot H^1(\R^2))$ for the same problem. In \cite{FK}
Fujita-Kato proved that for $u^0\in\dot H^{1/2}(\R^3)$ there exist
$T^*\in (0,+\infty]$ and at least one solution
$$
u\in{\mathcal C}([0,T^*),\dot H^{1/2}(\R^3))\;\;\; and\;\;\; t^{1/4}
\nabla u\in {\mathcal C}([0,T^*),L^2(\R^3)).
$$
Moreover, if $\|u^0\|_{\dot H^{1/2}(\R^3)}<c\nu$ we have
$u\in{\mathcal C}_b(\R_+,\dot H^{1/2})\cap L^2(\R_+,\dot
H^{3/2})$.\\ In what follows, we summarize some classical and useful
results.
\begin{th1}\label{t00}\cite{FK} Let $u^0\in H^{1/2}(\R^3)$ be a divergence-free vectors field. There
exists $T>0$ and a unique solution $u\in\mathcal C([0,T], H^{1/2}(\R^3))\cap
L^2([0,T], H^{3/2}(\R^3))$. Moreover, if $\|u^0\|_{\dot
H^{1/2}(\R^3)}\leq c\nu$ we have $u\in{\mathcal C}_b(\R^+,
H^{1/2})\cap L^2(\R^+,\dot H^{3/2})$.
\end{th1}
\begin{th1}\label{t01}\cite{KAT} If $u^0\in H^s(\R^3)$, with $s>5/2$,
a divergence-free vectors field, then there exists a time $T>0$ and
a strong solution $u$ of ($NS_\nu$)  in ${\mathcal
C}([0,T];H^{s})\cap {\mathcal C}^1([0,T];H^{s-2})$.
\end{th1}
With similar hypothesis, we have the blow-up result for the
maximal solutions. \begin{th1}\label{t02}\cite{BKM}\cite{KPO} Let
$s>5/2$ and $u\in {\mathcal C}([0,T^*);H^{s})$ a solution of
$(NS_\nu)$ with $u\notin{\mathcal C}([0,T^*];H^{s})$. Then
$$ \int_0^{T^*}\|\omega(t)\|_{L^\infty}dt=+\infty,$$
where $\omega={curl} u\simeq \frac{1}{2}(\nabla u-^t\nabla u)$.\\
\end{th1}
\noindent In this work, we investigate  some effects of the elliptic
operator $-\nu\Delta$ on the solution, we prove that if $u^0\in
H^{1/2}$, then $u\in\mathcal C((0,T],H^s)$, for all $s\in\R$, and we
present some asymptotic behavior near  0 and $+\infty$ for the
global solutions. The first result in this direction is due to J.- Y. Chemin in \cite{C4}. Precisely
\begin{th1}\cite{C4} Let $u^0\in H^{1/2}(\R^3)$ be a divergence-free vector field. There exist a time $T>0$ and a unique solution $u$ to $(NS_\nu)$ satisfying$$\forall 0\leq t\leq T,\;\int_0^t\int_{\R^3}|\xi|^2e^{(\nu \tau)^{1/2}|\xi|}|\hat{u}(\tau,\xi)|^2d\tau d\xi\leq 4\|\nabla e^{\tau\frac{\nu}{2}\Delta}u^0\|_{L^2([0,t]\times\R^3)}^2.$$
Moreover, a constant $c$ exist such that, if $u^0$  satisfy $\|u^0\|_{\dot H^{1/2}}<c\nu$, then
$$\int_{\R^+\times\R^3}|\xi|^2e^{(\nu t)^{1/2}|\xi|}|\hat{u}(t,\xi)|^2dtd\xi\leq \frac{4}{\nu}\|u^0\|_{L^2(\R^3)}^2.$$
\end{th1}
Our main results are the following
\begin{th1}\label{t1}(Small initial data) Let  $u^0\in H^{1/2}(\R^3)$
a divergence-free vector field, such that $\|u^0\|_{\dot
H^{1/2}}\leq c\nu$, then there exists a unique $u_\nu\in\mathcal
C_b(\R^+,H^{1/2}(\R^3))\cap L^2(\R^+,\dot H^{3/2}(\R^3))$. Moreover,
for all $T>0$, there exists $\varepsilon=\varepsilon(T,\nu,u^0)>0$,
such that

$$\forall 0\leq t\leq T,\;\int_{\R^3}|\xi|e^{2\varepsilon\nu t|\xi|}|\hat{u}(t,\xi)|^2d\xi+\nu\int_0^t\int_{\R^3}|\xi|^3e^{2\varepsilon\nu \tau|\xi|}|\hat{u}(\tau,\xi)|^2d\tau d\xi\leq 2\|u^0\|_{\dot H^{1/2}}^2.$$
In addition, we
have
$$\begin{array}{lll}&&(\alpha_1)\;\forall s\in\mathbb R,
\;\;u_\nu\in \mathcal C(\mathbb R^{+*},H^s),\\\\
&&(\alpha_2)\;\forall s>1/2,\;\;t^{s-\frac{1}{2}}\|u_\nu(t)\|_{H^s}\leq C_s,\;\;t\rightarrow 0^+,\\\\
&&(\alpha_3)\;\forall 0<q<1/8,\;\|u_\nu(t)\|_{\dot H^{1/2}}\leq C_q t^{-q},\;t\rightarrow+\infty,\\\\
&&(\alpha_4)\;\forall s>1/2,\;\|u_\nu(t)\|_{\dot H^s}\leq
C_st^{-(s-\frac{1}{2})/2},\;t\rightarrow+\infty.\\
\end{array}$$
\end{th1}
\begin{rem} The property $(\alpha_1)$ imply that
$\forall t>0,\;u_\nu(t,.)\in \mathcal C^\infty(\mathbb R^3)$.
\end{rem}
\begin{th1}\label{t2}(Large initial data)
For all vector field $u^0\in H^{1/2}(\R^3)$ divergence-free, there
exists $T>0$ and a unique $u_\nu\in\mathcal C([0,T],
H^{1/2}(\R^3))\cap L^2([0,T],\dot H^{3/2}(\R^3))$. Moreover, there
exists $\varepsilon=\varepsilon(T,\nu,u^0)>0$, such that $$\forall 0\leq t\leq T,\;\int_{\R^3}|\xi|e^{2\varepsilon\nu t|\xi|}|\hat{u}(t,\xi)|^2d\xi+\nu\int_0^t\int_{\R^3}|\xi|^3e^{2\varepsilon\nu \tau|\xi|}|\hat{u}(\tau,\xi)|^2d\tau d\xi\leq 2\|u^0\|_{\dot H^{1/2}}^2.$$
In addition, we have
$$\begin{array}{lcl}&&(\beta_1)\;\forall s\in\mathbb R,\;\;u_\nu\in
\mathcal C((0,T],H^s),\\\\
&&(\beta_2)\;\forall s>1/2,\;\;\|u_\nu(t)\|_{H^s}\leq \frac{C}{t^{s-\frac{1}{2}}},\;\;t\rightarrow 0^+,\\\\
&&(\beta_3)\;\forall 0<t\leq T,\;\;u_\nu(t,.)\in \mathcal C^\infty(\mathbb R^3).\\
\end{array}$$\\
\end{th1}
\noindent The rest of this paper is organized as follows. Section 2
contains some notations and definitions. In Section 3, we
prove Theorem \ref{t1}. The proof is based on Friedrich methods,
classical product law in the homogeneous Sobolev spaces, classical
compactness methods and  elementary technical results. Section 4 is
devoted to prove Theorem \ref{t2}. The proof is inspired from  the
previews  section. In  last Section we derive some general proprieties of strong solutions of Navier-Stokes equations
$(NS_\nu)$.
\section{Notations}
\noindent In this short section we collect some notations and
definitions that will be used later on.\\
$\bullet$ The Fourier transformation is normalized as
$${\mathcal
F}(f)(\xi)=\stackrel{\wedge}{f}(\xi)=\displaystyle\int_{\R^3}\exp(-ix.\xi)f(x)dx,\;\;\;
\xi=(\xi_1,\xi_2,\xi_3)\in\R^3.$$ $\bullet$ The inverse Fourier
formula is
$${\mathcal
F}^{-1}(g)(x)=\displaystyle(2\pi)^{-3}\int_{\R^3}\exp(i\xi.x)f(\xi)d\xi,\;\;\;
x=(x_1,x_2,x_3)\in\R^3.$$ $\bullet$ For $s\in\R$, $H^s(\R^3)$
denotes the usual non homogeneous Sobolev space on $\R^3$ and
$<.,.>_{H^s(\R^3)}$ denotes the usual scalar product on
$H^s(\R^3)$.\\\\
$\bullet$ For $s\in\R$, $\dot H^s(\R^3)$ denotes the usual
homogeneous Sobolev space on $\R^3$ and $<.,.>_{\dot H^s(\R^3)}$
denotes the usual scalar product on
$\dot H^s(\R^3)$.\\\\
$\bullet$ The convolution product of a suitable pair of functions
$f$ and $g$ on $\R^3$ is given by $$(f*g)(x):=\int_{\mathbb R^3}f(y)g(x-y)dy.$$\\
$\bullet$ For any Banach space $(B,\|.\|)$, any real number $1\leq
p\leq \infty$ and any time $T>0$, we will denote by $L^p_T(B)$ the
space of all measurable functions
$$t\in[0,T]\rightarrow f(t)\in
B$$such that $$(t\rightarrow \|f(t)\|)\in L^p([0,T]).$$\\
$\bullet$ If $f=(f_1,f_2,f_3)$ and $g=(g_1,g_2,g_3)$ are two vector
fields, we set $$f\otimes g:=(g_1f,g_2f,g_3f),$$ and$$div(f\otimes
g):=(div(g_1f),div(g_2f),div(g_3f)).$$\\
$\bullet$ For any subset $X$ of a set $E$, the symbol ${\bf 1}_X$
denote the characteristic function of $X$ defined by $${\bf
1}_X(x)=1\;\; \mbox{if}\;\; x\in X,\;\;\;{\bf 1}_X(x)=0\;
\mbox{elsewhere}.$$
\section{Proof of Theorem \ref{t1}}
\noindent We begin by recalling a fundamental lemma concerning some product laws in Sobolev spaces.
\begin{lem}\label{l1}(see \cite{C0}) Let $s,s'$ tow reals
numbers such that
$$s<3/2\qquad \mbox{ and } \qquad  s+s'>0.$$There exists a positive constant $C:=C(s,s')$, such that for
all $f,g\in\dot H^s(\mathbb R^3)\cap\dot H^{s'}(\mathbb R^3)$,
$$\|fg\|_{\dot H^{s+s'-\frac{3}{2}}(\mathbb R^3)}\leq C
\Big(\|f\|_{\dot H^{s}(\mathbb R^3)}\|g\|_{\dot H^{s'}(\mathbb
R^3)}+ \|f\|_{\dot H^{s'}(\mathbb R^3)}\|g\|_{\dot H^{s}(\mathbb
R^3)} \Big).$$ If $s,s'<3/2$ and $s+s'>0$, there exist a constant
$c=c(s,s')$,
$$\|fg\|_{\dot H^{s+s'-\frac{3}{2}}(\mathbb R^3)}\leq c
\|f\|_{\dot H^{s}(\mathbb R^3)}\|g\|_{\dot H^{s'}(\mathbb
R^3)}.$$\\
\end{lem}
\noindent For a strictly positive integer $n$, the Friedrich's
operator $J_n$ is defined by
$$J_n(f):={\mathcal F}^{-1}\Big({\bf
1}_{\{|\xi|<n\}}{\mathcal F}(f)\Big).$$ Consider the following
approximate Navier-Stokes system $(NS_{n,\nu})$ on
$\R_{+}\times\R^3$,
$$ \left\{\begin{array}{l}
  \displaystyle\partial_t
u-  \nu\Delta J_n u+J_n div(J_nu\otimes J_n u)=\nabla\Delta^{-1}J_n
div(J_nu\otimes J_n u),\\\\
u_{\mid t=0}=J_nu_0.\\
\end{array}\right.$$ Then by the ordinary differential
equations theory  the system $(NS_{n,\nu})$ has a unique maximal
solution $u_{n,\nu}$ in the space ${\mathcal
C}^1([0,T_n^*),L^2(\R^3)).$ Using the uniqueness and the fact
$J_n^2=J_n$ we obtain
$$\left\{\begin{array}{lcl}
J_nu_{n,\nu}&=&u_{n,\nu}\\\\
div(u_{n,\nu}(t))&=&0,\;\; \forall t\in [0,T_n^*),
\end{array}\right.
$$
hence $u_{n,\nu}$ satisfies
$$ \left\{\begin{array}{l}
  \displaystyle\partial_t
u_{n,\nu}-  \nu\Delta u_{n,\nu}+J_n div (u_{n,\nu}\otimes
u_{n,\nu})=\nabla\Delta^{-1}J_n
div(u_{n,\nu}\otimes u),\\\\
{u_{n,\nu}}_{\mid t=0}=J_nu_0.\\
\end{array}\right.$$
Taking the scalar product in $L^2(\R^3)$, we obtain, for $t\in
[0,T_n^*)$,
$$
\partial_t \|u_{n,\nu}\|_{L^2 }^2+2\nu\|\nabla
u_{n,\nu}\|_{ L^2}^2\leq 0,
$$
then it follows for all $t\in [0,T_n^*)$,
$$\|u_{n,\nu}(t)\|_{L^2}^2\leq \|u^0\|_{L^2}^2,$$ which implies
$T_n^*=+\infty$,  and the following estimate holds  for all $t\geq
0$
\begin{equation}\label{es21} \|u_{n,\nu}(t)\|_{L^2}^2+2\nu\|\nabla
u_{n,\nu}\|_{L^2_t(L^2)}^2\leq \|u^0\|_{L^2}^2.
\end{equation}
Let $T>0$ a fixed time, and let
$$(C_1)\;\;\;\;\;\;\;\;\;\;\varepsilon:=\min\Big(\frac{1}{2},c\nu,(c\nu)^2,CT^{-1/2},CT^{-1/3}\Big),$$
with $C$ depending only of $\nu,\; u^0$.\\\\ Using Fourier
transformation we obtain from the above system \\
$$(E_1):\, \partial_t\stackrel{\wedge} {v}_{n,\nu,\varepsilon}
+\nu|\xi|(|\xi|-\varepsilon)\stackrel{\wedge}
{v}_{n,\nu,\varepsilon}+e^{\varepsilon\nu t|\xi|}\mathcal
F\Big(J_n(u_{n,\nu}.\nabla u_{n,\nu})\Big)=e^{\varepsilon\nu
t|\xi|}\mathcal F\Big(\nabla\Delta^{-1}J_n (u_{n,\nu}.\nabla
u_{n,\nu})\Big),$$\\ where $v_{n,\nu,\varepsilon}:=\mathcal
F^{-1}\Big(e^{\varepsilon\nu t|\xi|}\stackrel{\wedge}
{u}_{n,\nu}\Big)$. Take $\overline{E}_1.\stackrel{\wedge}
{v}_{n,\nu,\varepsilon}+E_1.\overline{\stackrel{\wedge}
{v}}_{n,\nu,\varepsilon}$, we obtain
$$\partial_t|\stackrel{\wedge}
{v}_{n,\nu,\varepsilon}|^2
+2\nu|\xi|(|\xi|-\varepsilon)|\stackrel{\wedge}
{v}_{n,\nu,\varepsilon}|^2=-2\Re e\Big(e^{\varepsilon\nu
t|\xi|}\mathcal F\Big(J_n (u_{n,\nu}.\nabla
u_{n,\nu})\Big)\overline{\stackrel{\wedge}
{v}}_{n,\nu,\varepsilon}\Big).$$ Using the following elementary inequality
$$e^{a|\xi|}\leq e^{a|\xi-\eta|}e^{a|\eta|},\;\;\forall
a\in\mathbb R^+,\;\;\forall \xi,\eta\in\mathbb R^3,$$ we obtain
$$(E_2): \, \partial_t|\stackrel{\wedge} {v}_{n,\nu,\varepsilon}|^2
+2\nu|\xi|(|\xi|-\varepsilon)|\stackrel{\wedge}
{v}_{n,\nu,\varepsilon}|^2\leq 2|\stackrel{\wedge}
{v}_{n,\nu,\varepsilon}|\ast|\stackrel{\wedge}
{v}_{n,\nu,\varepsilon}|.|\stackrel{\wedge} {\nabla
v}_{n,\nu,\varepsilon}|.$$ To make good estimates, we decompose
$v_{n,\nu,\varepsilon}$ as
$v_{n,\nu,\varepsilon}=X_{n,\nu,\varepsilon}+Y_{n,\nu,\varepsilon}$,
with $$\begin{array}{lcl}  {X}_{n,\nu,\varepsilon}&=&\mathcal
F^{-1}\Big({\bf 1}_{\{|\xi|>2\varepsilon\}}\stackrel{\wedge}
{v}_{n,\nu,\varepsilon}\Big)\\
{Y}_{n,\nu,\varepsilon}&=&\mathcal F^{-1}\Big({\bf
1}_{\{|\xi|\leq 2\varepsilon\}}\stackrel{\wedge}
{v}_{n,\nu,\varepsilon}\Big).\\\\
\end{array}$$
{\bf Estimate of $\|X_{n,\nu,\varepsilon}\|_{\dot H^{1/2}}^2$:} Multiply $(E_2)$ by $|\xi|$ and integrate over
$\{|\xi|>2\varepsilon\}$, we obtain
$$\partial_t\|X_{n,\nu,\varepsilon}\|_{\dot H^{1/2}}^2
+\nu\|\nabla X_{n,\nu,\varepsilon}\|_{\dot H^{1/2}}^2\leq
\||\xi|^{1/2}|\stackrel{\wedge}
{v}_{n,\nu,\varepsilon}|\ast|\stackrel{\wedge}
{v}_{n,\nu,\varepsilon}|\|_{L^2}\|\nabla
X_{n,\nu,\varepsilon}\|_{\dot H^{1/2}}\leq\sum_{j=1}^3I_j,$$ where
$$\begin{array}{lcl}
I_1&=&\||\xi|^{1/2}|\stackrel{\wedge}
{X}_{n,\nu,\varepsilon}|\ast|\stackrel{\wedge}
{X}_{n,\nu,\varepsilon}|\|_{L^2}\|\nabla
X_{n,\nu,\varepsilon}\|_{\dot H^{1/2}},\\
I_2&=&2\||\xi|^{1/2}|\stackrel{\wedge}
{X}_{n,\nu,\varepsilon}|\ast|\stackrel{\wedge}
{Y}_{n,\nu,\varepsilon}|\|_{L^2}\|\nabla
X_{n,\nu,\varepsilon}\|_{\dot H^{1/2}},\\
I_3&=&\||\xi|^{1/2}|\stackrel{\wedge}
{Y}_{n,\nu,\varepsilon}|\ast|\stackrel{\wedge}
{Y}_{n,\nu,\varepsilon}|\|_{L^2}\|\nabla
X_{n,\nu,\varepsilon}\|_{\dot H^{1/2}}.\\\\ \end{array}$$ Using now
Lemma \ref{l1}, we get
$$
I_1\leq C\|X_{n,\nu,\varepsilon}\|_{\dot H^{1/2}} \|\nabla
X_{n,\nu,\varepsilon}\|_{\dot H^{1/2}}^2.$$ Once again, Lemma
\ref{l1}, combined together with  the energy estimate (\ref{es21}),
give
$$\begin{array}{lcl}I_2&\leq&C\|X_{n,\nu,\varepsilon}\|_{\dot
H^{1/2}}\|\nabla Y_{n,\nu,\varepsilon}\|_{\dot H^{1/2}}\|\nabla
X_{n,\nu,\varepsilon}\|_{\dot H^{1/2}}
+C\|Y_{n,\nu,\varepsilon}\|_{\dot H^{1/2}}
\|\nabla X_{n,\nu,\varepsilon}\|_{\dot H^{1/2}}^2\\\\
&\leq&
C\varepsilon^{3/2}\|Y_{n,\nu,\varepsilon}\|_{L^2}\|X_{n,\nu,\varepsilon}\|_{\dot
H^{1/2}}\|\nabla X_{n,\nu,\varepsilon}\|_{\dot
H^{1/2}}+C\varepsilon^{1/2}\|Y_{n,\nu,\varepsilon}\|_{L^2}
\|\nabla X_{n,\nu,\varepsilon}\|_{\dot H^{1/2}}^2\\\\
&\leq&
C\varepsilon^{3/2}\|u^0\|_{L^2}\|X_{n,\nu,\varepsilon}\|_{\dot
H^{1/2}}\|\nabla X_{n,\nu,\varepsilon}\|_{\dot
H^{1/2}}+C\varepsilon^{1/2}\|u^0\|_{L^2}
\|\nabla X_{n,\nu,\varepsilon}\|_{\dot H^{1/2}}^2\\\\
&\leq& C\varepsilon^3\|u^0\|_{L^2}^2\|X_{n,\nu,\varepsilon}\|_{\dot
H^{1/2}}^2+\displaystyle(\frac{\nu}{20}+C\varepsilon^{1/2}\|u^0\|_{L^2}
)\|\nabla
X_{n,\nu,\varepsilon}\|_{\dot H^{1/2}}^2\\\\
&\leq& \displaystyle\frac{1}{10}\|X_{n,\nu,\varepsilon}\|_{\dot
H^{1/2}}^2+\displaystyle\frac{\nu}{10}\|\nabla
X_{n,\nu,\varepsilon}\|_{\dot H^{1/2}}^2.\\\\
\end{array}$$
Similary, we obtain
$$\begin{array}{lcl}
I_3&\leq&C\varepsilon^2\|Y_{n,\nu,\varepsilon}\|_{L^2}^2\|\nabla
X_{n,\nu,\varepsilon}\|_{\dot H^{1/2}}\\\\
&\leq&C\varepsilon^2\|u^0\|_{L^2}^2\|\nabla
X_{n,\nu,\varepsilon}\|_{\dot H^{1/2}}\\\\
&\leq&C\varepsilon^4\|u^0\|_{L^2}^4+\displaystyle\frac{\nu}{10}\|\nabla
X_{n,\nu,\varepsilon}\|_{\dot H^{1/2}}^2\\\\
&\leq&\displaystyle\frac{\|u^0\|_{\dot
H^{1/2}}^2}{10T}+\displaystyle\frac{\nu}{10}\|\nabla
X_{n,\nu,\varepsilon}\|_{\dot H^{1/2}}^2,\\\\
\end{array}$$
which leads to
$$\begin{array}{lcl}\partial_t\|X_{n,\nu,\varepsilon}\|_{\dot
H^{1/2}}^2 +\nu\|\nabla X_{n,\nu,\varepsilon}\|_{\dot
H^{1/2}}^2&\leq&C\|X_{n,\nu,\varepsilon}\|_{\dot H^{1/2}}\|\nabla
X_{n,\nu,\varepsilon}\|_{\dot H^{1/2}}^2\\\\
&+&C_2\varepsilon^3\nu^{-1}\|X_{n,\nu,\varepsilon}\|_{\dot
H^{1/2}}^2+(C_2\varepsilon^{1/2}+\frac{\nu}{10})\|\nabla
X_{n,\nu,\varepsilon}\|_{\dot H^{1/2}}^2\\\\
&+&C_2\varepsilon^2,\\\\
\end{array}$$
and hence $$\partial_t\|X_{n,\nu,\varepsilon}\|_{\dot H^{1/2}}^2
+\frac{\nu}{2}\|\nabla X_{n,\nu,\varepsilon}\|_{\dot H^{1/2}}^2\leq
C\|X_{n,\nu,\varepsilon}\|_{\dot H^{1/2}}\|\nabla
X_{n,\nu,\varepsilon}\|_{\dot
H^{1/2}}^2+C_2\varepsilon^3\nu^{-1}\|X_{n,\nu,\varepsilon}\|_{\dot
H^{1/2}}^2+C_2\varepsilon^2.$$ Let a time $T_{n,\nu,\varepsilon}$
define by
$$T_{n,\nu,\varepsilon}:=\sup\{t\geq
0,\;\|X_{n,\nu,\varepsilon}\|_{L^\infty_t(\dot H^{1/2})}<2c\nu\}.$$
For $0\leq t<\min(T,T_{n,\nu,\varepsilon})$, we have
$$\|X_{n,\nu,\varepsilon}(t)\|_{\dot H^{1/2}}^2
+\frac{\nu}{4}\int_0^t\|\nabla X_{n,\nu,\varepsilon}\|_{\dot
H^{1/2}}^2\leq \|X_{n,\nu,\varepsilon}(0)\|_{\dot H^{1/2}}^2+
C\varepsilon^2t+
C\varepsilon^3\nu^{-1}\int_0^t\|X_{n,\nu,\varepsilon}\|_{\dot
H^{1/2}}^2,$$ and by Gronwall's lemma we get
$$\begin{array}{lcl}\|X_{n,\nu,\varepsilon}\|_{L^\infty_t(\dot H^{1/2})}^2
&\leq&(\|u^0\|_{\dot H^{1/2}}^2+
C\varepsilon^2T)e^{C\varepsilon^3\nu^{-1}T}\\\\
&\leq&\displaystyle\frac{3}{2}\|u^0\|_{\dot H^{1/2}}^2<(2c\nu)^2,
\end{array}$$
that is  $T_{n,\nu,\varepsilon}>T$, and for all $0\leq t\leq T$,
\begin{equation}\label{exx}\|X_{n,\nu,\varepsilon}\|_{L^\infty_t(\dot
H^{1/2})}^2+\frac{\nu}{4}\|\nabla
X_{n,\nu,\varepsilon}\|_{L^2_t(\dot H^{1/2})}^2\leq
\displaystyle\frac{3}{2}\|u^0\|_{\dot H^{1/2}}^2.\end{equation} {\bf
Estimate of $\| Y_{n,\nu,\varepsilon}\|_{\dot H^{1/2}}^2$:} To
estimate $\| Y_{n,\nu,\varepsilon}\|_{\dot H^{1/2}}^2$ we integrate
$(E_2)$ over $\{|\xi|<2\varepsilon\}$ to obtain
$$\begin{array}{lcl}\partial_t\|Y_{n,\nu,\varepsilon}\|_{\dot H^{1/2}}^2
+\nu\|\nabla {Y}_{n,\nu,\varepsilon}\|_{\dot H^{1/2}}^2&\leq&
2\||\xi|^{1/2}|\stackrel{\wedge} {u}_{n,\nu}|\ast|\stackrel{\wedge}
{u}_{n,\nu,\varepsilon}|\|_{L^2(B(0,2\varepsilon))}.
\|\nabla {Y}_{n,\nu,\varepsilon}\|_{\dot H^{1/2}}\\
&\leq&C\displaystyle\Big(\int_{B(0,2\varepsilon)}|\xi|\Big)^{1/2}d\xi\||\stackrel{\wedge}
{u}_{n,\nu}|\ast|\stackrel{\wedge}
{u}_{n,\nu}|\|_{L^\infty}.\|\nabla {Y}_{n,\nu}\|_{\dot H^{1/2}},\\\\
\end{array}$$
and by Young inequality it follows
$$\begin{array}{lcl}\partial_t\|Y_{n,\nu,\varepsilon}\|_{\dot H^{1/2}}^2
+2\nu\|\nabla Y_{n,\nu,\varepsilon}\|_{\dot H^{1/2}}^2&\leq&
C\varepsilon^2\|\stackrel{\wedge} {u}_{n,\nu}\|_{L^2}^2\|\nabla
Y_{n,\nu}\|_{\dot H^{1/2}}\\\\
&\leq& C\|u^0\|_{L^2}^2\varepsilon^2\|\nabla Y_{n,\nu,\varepsilon}\|_{\dot H^{1/2}}\\\\
&\leq&C\nu^{-1}\|u^0\|_{L^2}^4\varepsilon^4+\nu\|\nabla
Y_{n,\nu,\varepsilon}\|_{\dot H^{1/2}}^2.
\end{array}$$
An easy computation shows that  for all $t\in [0,T]$
\begin{eqnarray}\label{exy}\nonumber
\|Y_{n,\nu,\varepsilon}(t)\|_{\dot H^{1/2}}^2+\nu\|\nabla
Y_{n,\nu,\varepsilon}\|_{L_t^\infty(\dot H^{1/2})}^2&\leq&
\|Y_{n,\nu,\varepsilon}(0)\|_{\dot
H^{1/2}}^2+C\varepsilon^4t\\
&\leq&\nonumber\varepsilon\|u_o\|_{L^2}^2+C\varepsilon^4T\\
&\leq&\displaystyle\frac{1}{4}\|u_o\|_{L^2}^2.
\end{eqnarray}
Thanks to equations  (\ref{exx})-(\ref{exy})  we obtain  for all
$t\in[0,T]$,
$$\|v_{n,\nu,\varepsilon}(t)\|_{\dot
H^{1/2}}^2+\nu\displaystyle\int_0^t\|\nabla
v_{n,\nu,\varepsilon}\|_{\dot H^{1/2}}^2 \leq 2\|u^0\|_{\dot
H^{1/2}}^2 .$$ Finally, a standard compactness argument gives the
global existence result, precisely : There exists $u_\nu\in \mathcal
C_b(\R^+,H^{1/2}(\mathbb R^3))\cap L^2(\R^+,\dot H^{3/2}(\mathbb
R^3))$ of $(NS_\nu)$, such that for all $T>0$, and
$\varepsilon=\varepsilon_{_T}:=\min(1/2,c\nu,(c\nu)^2,CT^{-1/2},CT^{-1/3})$
we have \begin{equation}\label{exp1}\forall
t\in[0,T],\;\int_{\xi}e^{2\varepsilon\nu
t|\xi|}|\xi|.|\stackrel{\wedge}
{u}_\nu(t,\xi)|^2d\xi+\nu\int_0^t\int_{\xi}e^{2\varepsilon\nu
\tau|\xi|}|\xi|^3|\stackrel{\wedge} {u}_\nu(\tau,\xi)|^2d\xi
d\tau\leq 2\|u^0\|_{\dot H^{1/2}}^2.\end{equation} The equation
(\ref{es21}) yields
 \begin{equation}\label{exp2}\forall
t\geq
0,\;\;\;\;\|u_\nu(t)\|_{L^2}^2+2\nu\displaystyle\int_0^t\|\nabla
u_\nu\|_{L^2}^2 \leq \|u^0\|_{L^2}^2,\end{equation} which implies
$(\alpha_1)-(\alpha_2)$.\\\\
Proof of $(\alpha_3)$. For $a>0$ we have
$$\begin{array}{lcl}\|u_\nu(T)\|_{\dot H^{1/2}}&\leq& \displaystyle
\Big(\int_{\{|\xi|<a\}}|\xi||\stackrel{\wedge}{u}_\nu(T,\xi)|^2\Big)^{1/2}
+\Big(\int_{\{|\xi|>a\}}|\xi||\stackrel{\wedge}{u}_\nu(T,\xi)|^2\Big)^{1/2}\\\\
&\leq&\displaystyle \sqrt{a}\|u_\nu(t)\|_{L^2}
+e^{-\varepsilon_{_T}\nu Ta^2}\Big(\int_{\{|\xi|>a\}}|\xi|e^{\varepsilon_{_T}\nu T|\xi|^2}|\stackrel{\wedge}{u}_\nu|^2\Big)^{1/2}\\\\
&\leq&\displaystyle \sqrt{a}\|u^0\|_{L^2} +e^{-\varepsilon_{_T}\nu
Ta^2}\sqrt{3}\|u^0\|_{\dot H^{1/2}}.
\end{array}$$
For  large naught  time $T$, we have $\varepsilon_{_T}=CT^{-1/2}$ and choosing  $a=T^{-r}$ with $0<r<\frac{1}{4}$, it follows
$$\|u_\nu(T)\|_{\dot H^{1/2}}\leq \frac{\|u^0\|_{L^2}}{T^{r/2}}
+e^{-CT^{\frac{1-4r}{2}}}\sqrt{3}\|u^0\|_{\dot H^{1/2}}.$$ as desired.\\\\
The case of  $(\alpha_4)$, follows as well by choosing
$a:=\frac{C}{T^{1/2}}$. The proof of theorem $\ref{t1}$ is
completed.
\endproof
\section{Proof of theorem \ref{t2}}\noindent We apply the Friedrich's method's, we obtain the
existence and uniqueness of solution $u_{n,\nu}\in{\mathcal
C}^1(\mathbb R^+,L^2(\R^3))$ of the following system
$$ \left\{\begin{array}{l}
  \displaystyle\partial_t
u-  \nu\Delta u+J_n div (u\otimes u)=\nabla\Delta^{-1}J_n
div(u\otimes u),\\\\
u_{\mid t=0}=J_nu_0.\\
\end{array}\right.$$
And we have
\begin{equation}
\forall t\geq 0\;\;\;\|u_{n,\nu}(t)\|_{L^2}^2+2\nu\|\nabla
u_{n,\nu}\|_{L^2_t(L^2)}^2\leq \|u^0\|_{L^2}^2.
\end{equation}
For simplification we don't well noting the index $n$.\\\\
Let $N\in\mathbb N$, such that
$$\|\mathcal F^{-1}({\bf 1}_{\{|\xi|<2^{-N}\}\cup
\{|\xi|>2^N\}}\mathcal F(u^0))\|_{H^{1/2}}<\min(c\nu,c\nu^{3/2}).$$
And define $T$, $u^0_N$, $v_{N,L}$, $w_{N,\nu}$, by
$$(C_2)\;\;\;\;T=T(\nu,u^0):=-\nu^{-1}2^{-N}\log\Big(1-\min(1/2,\nu\|u^0\|_{\dot
H^{1/2}}^{-4}\min(c\nu,(c\nu)^{2/3}))\Big)>0.$$
$$\begin{array}{lcl}u^0_N:&=&\mathcal F^{-1}({\bf
1}_{\{2^{-N}<|\xi|<2^N\}}\mathcal F(u^0))\\\\
v_{N,L}:&=&e^{\nu t\Delta}u^0_N\\\\
w_{N,\nu}:&=&u_{n,\nu}-v_{N,L}.
\end{array}$$
We have
$$\partial_tw_{N,\nu}-\nu\Delta w_{N,\nu}+w_{N,\nu}.\nabla
w_{N,\nu}+v_{N,L}.\nabla w_{n,\nu}+w_{n,\nu}.\nabla
v_{N,L}=-v_{N,L}.\nabla v_{N,L}.$$ We introduce the real number
$\varepsilon$ defined by
$$(C_3)\;\;\;\;\;\;\quad\quad\varepsilon:=\min\Big(\frac{1}{2.2^N},(10CT)^{-1/3},(\min(c\nu,(c\nu)^{3/2}))^{1/4}(10CT)^{-1/4}\Big)>0,$$
and
$$\begin{array}{lcl}
V_{N,L}&=&e^{\varepsilon\nu t|D|}v_{N,L}\\\\
U_{\nu}&=&e^{\varepsilon\nu t|D|}u_{\nu}\\\\
W_{N,\nu}&=&e^{\varepsilon\nu t|D|}w_{N,\nu}\\\\
\alpha_{N,\nu,\varepsilon}&=&{\bf 1}_{\{|D|>2\varepsilon\}}W_{N,\nu}\\\\
\beta_{N,\nu,\varepsilon}&=&{\bf 1}_{\{|D|\leq
2\varepsilon\}}W_{N,\nu}.
\end{array}$$
Arguing as in the last section, we obtain
$$\partial_t\|\alpha_{N,\nu,\varepsilon}\|_{\dot H^{1/2}}^2
+\nu\|\nabla \alpha_{N,\nu,\varepsilon}\|_{\dot H^{1/2}}^2\leq
\sum_{j=1}^3K_j,$$ with
$$\begin{array}{lcl}
K_1&=&\displaystyle\int_\xi|\xi||\stackrel{\wedge}
{W}_{N,\nu}|\ast|\stackrel{\wedge}
{W}_{N,\nu}|.|\stackrel{\wedge}{\nabla
\alpha}_{n,\nu,\varepsilon}|,\\\\
K_2&=&\displaystyle\int_\xi|\xi|^{-1/2}\Big(|\stackrel{\wedge}
{V}_{N,L}|\ast|\stackrel{\wedge} {\nabla
W}_{N,\nu}|+|\stackrel{\wedge} {\nabla
V}_{N,L}|\ast|\stackrel{\wedge}
{W}_{N,\nu}|\Big).|\xi|^{3/2}|\stackrel{\wedge}{
\alpha}_{n,\nu,\varepsilon}|,\\\\
K_3&=&\displaystyle\int_\xi|\xi|^{-1/2}.|\stackrel{\wedge}
{V}_{N,L}|\ast|\stackrel{\wedge} {\nabla
V}_{N,L}|.|\xi|^{3/2}|\stackrel{\wedge}{
\alpha}_{n,\nu,\varepsilon}|.
\end{array}$$
Using the Lemma \ref{l1}, we obtain
$$\begin{array}{lcl}
K_1&\leq&C\|W_{N,\nu}\|_{\dot H^{1/2}}\|\nabla W_{N,\nu}\|_{\dot
H^{1/2}}\|\nabla \alpha_{n,\nu,\varepsilon}\|_{\dot H^{1/2}},\\\\
K_2&\leq&C\|\nabla V_{N,L}\|_{L^{2}}\|\nabla W_{N,\nu}\|_{L^2}
\|\nabla \alpha_{n,\nu,\varepsilon}\|_{\dot H^{1/2}},\\\\
K_3&\leq&C\|\nabla V_{N,L}\|_{L^{2}}^2\|\nabla
\alpha_{n,\nu,\varepsilon}\|_{\dot H^{1/2}}.
\end{array}$$
{\bf Estimate of $K_1$:} We have
$$\begin{array}{lcl}
\|W_{N,\nu}\|_{\dot H^{1/2}}&\leq&\|\alpha_{N,\nu}\|_{\dot
H^{1/2}}+\|\beta_{N,\nu}\|_{\dot H^{1/2}}\leq
\|\alpha_{N,\nu}\|_{\dot H^{1/2}}+C\varepsilon^{1/2}\\\\
\|\nabla W_{N,\nu}\|_{\dot
H^{1/2}}&\leq&\|\nabla\alpha_{N,\nu}\|_{\dot
H^{1/2}}+\|\nabla\beta_{N,\nu}\|_{\dot H^{1/2}}\leq \|\nabla
\alpha_{N,\nu}\|_{\dot H^{1/2}}+C\varepsilon^{3/2},
\end{array}$$
then \begin{equation}\label{k1}K_1\leq
C\varepsilon^4+C_0\varepsilon^3\|\alpha_{N,\nu}\|_{\dot
H^{1/2}}^2+\Big(C_0\varepsilon^{1/2}+\frac{\nu}{10}+C\|\alpha_{N,\nu}\|_{\dot
H^{1/2}}\Big)\|\nabla\alpha_{N,\nu}\|_{\dot
H^{1/2}}^2.\end{equation} {\bf Estimate of $K_2$:} Similary, we
obtain
\begin{equation}\label{k2}K_2\leq \|\nabla V_{N,L}\|_{\dot
H^{1/2}}^4+C\varepsilon^4+C\varepsilon^3\|\alpha_{N,\nu}\|_{\dot
H^{1/2}}^2+\Big(C\varepsilon^{1/3}+\frac{\nu}{10}+C\|\alpha_{N,\nu}\|_{\dot
H^{1/2}}^{2/3}\Big)\|\nabla\alpha_{N,\nu}\|_{\dot
H^{1/2}}^2.\end{equation} {\bf Estimate of $K_3$:} We have
\begin{equation}\label{k3}K_3\leq C\|\nabla V_{N,L}\|_{L^2}^4
+\frac{\nu}{10}\|\nabla\alpha_{N,\nu}\|_{\dot
H^{1/2}}^2.\\\\\end{equation} Let$$T^*=\sup\{t>0,\;\;\|\alpha_{N,\nu}\|_{L_t^\infty(\dot
H^{1/2})}<2\min(c\nu,c\nu^{3/2})\}.$$ For $0\leq t<\min(T,T^*)$, we
have\\
$$\begin{array}{lcl}\|\alpha_{N,\nu,\varepsilon}\|_{L_t^\infty(\dot
H^{1/2})}^2 +\frac{\nu}{2}\|\nabla
\alpha_{N,\nu,\varepsilon}\|_{L_t^2(\dot H^{1/2})}^2&\leq&
\|\alpha_{N,\nu}(0)\|_{\dot H^{1/2}}^2+
C\varepsilon^4T+C\|\nabla V_{N,L}\|_{L_T^4(L^2)}^4\\\\
&+&C\varepsilon^3T\|\alpha_{N,\nu,\varepsilon}\|_{L_t^\infty(\dot
H^{1/2})}^2.\\\\
\end{array}$$ By inequalities (\ref{k1}), (\ref{k2}), (\ref{k3}) and the choices $(C_2)$, $(C_3)$, we obtain
$$\forall t\in (0,\min(T,T^*)),\;\;\|\alpha_{N,\nu,\varepsilon}\|_{L_t^\infty(\dot
H^{1/2})}^2 +\nu\|\nabla \alpha_{N,\nu,\varepsilon}\|_{L_t^2(\dot
H^{1/2})}^2\leq 2\min(c\nu,c\nu^{3/2})^2.$$ Then $T^*>T$, in
particular
$$\|\alpha_{N,\nu,\varepsilon}\|_{L_T^\infty(\dot
H^{1/2})}^2 +\nu\|\nabla \alpha_{N,\nu,\varepsilon}\|_{L_T^2(\dot
H^{1/2})}^2\leq 2\min(c\nu,c\nu^{3/2})^2,$$
$$\|W_{N,\nu}\|_{L_T^\infty(\dot H^{1/2})}^2
+\nu\|\nabla W_{N,\nu}\|_{L_T^2(\dot H^{1/2})}^2\leq
3\min(c\nu,c\nu^{3/2})^2,$$ and we can deduce
$$\|U_\nu\|_{L_T^\infty(\dot H^{1/2})}^2
+\nu\|\nabla U_\nu\|_{L_T^2(\dot H^{1/2})}^2\leq 2\|u^0\|_{\dot
H^{1/2}}^2.$$ Finally, a standard compactness argument gives the
local existence result. Moreover the solution satisfies
$(\beta_1)-(\beta_2)$. This achieved the proof of Theorem \ref{t2}.
\endproof
\section{General Proprieties Of Strong Solutions}
\noindent This section combines the previous results and Theorems
\ref{t01}, \ref{t02} to derive some proprieties of any strong
solutions of Navier-Stokes equations. The precise statements are the following. \begin{th1}\label{apt1} If $u\in \mathcal C([0,T_0],H^{1/2})\cap L^2([0,T_0],H^{3/2})$ is a
solution of $(NS_\nu)$, then\\
$$\begin{array}{lcl}&&(\beta_1')\;\forall s\in\mathbb R,\;\;u\in
\mathcal C(]0,T_0],H^s),\\\\
&&(\beta_2')\;\forall s>1/2,\;\;t^{s-\frac{1}{2}}\|u(t)\|_{H^s}\leq C_s,\;\;t\rightarrow 0^+,\\\\
&&(\beta_3')\;\forall t\in ]0,T_0],\;\;u(t,.)\in \mathcal C^\infty(\mathbb R^3).\\\\
\end{array}$$
\end{th1}
\begin{th1}\label{apt2} If $u\in \mathcal C(\R^+,H^{1/2})\cap L^2(\R^+,\dot H^{3/2})$ is a
solution of $(NS_\nu)$, then $u\in L^\infty (\R^+,H^{1/2})$, and\\
$$\begin{array}{lcl}&&(\alpha_1')\;\forall s\in\mathbb R,\;\;u\in
\mathcal C(\R^{*+},H^s),\\\\
&&(\alpha_2')\;\forall s>1/2,\;\;t^{s-\frac{1}{2}}\|u(t)\|_{H^s}\leq C_s,\;\;t\rightarrow 0^+,\\\\
&&(\alpha_3')\;\forall 0<q<1/8,\;\|u(t)\|_{\dot H^{1/2}}\leq C_q t^{-q},\;t\rightarrow+\infty,\\\\
&&(\alpha_4')\;\forall s>1/2,\;\|u(t)\|_{\dot H^s}\leq
C_st^{-(s-\frac{1}{2})/2},\;t\rightarrow+\infty.\\\\
\end{array}$$
\end{th1}
\begin{rem} $(\beta_3')$ is an easy consequence of $(\beta_1')$.\\
\end{rem}
\subsection{Proof Of Theorem \ref{apt1}}
\noindent Using Theorem \ref{t1}, there exists $T>0$ (suppose that
$T<T_0$), and $v_1\in\mathcal C([0,T],H^{1/2})\cap L^2_T(H^{3/2})$,
satisfying $(\beta_1)$. By the uniqueness, we have $v_1=u$ on
$[0,T]$. Let $s>5/2$, and we consider the following system
$$ \left\{\begin{array}{l}
  \displaystyle\partial_t
v-  \nu\Delta v+ v.\nabla v=-\nabla p,\quad
\mbox{on}\;\; \R_+\times\R^3, \\\\
div\; (v) = 0 \quad
\mbox{on}\;\; \R_+\times\R^3,  \\\\
v_{\mid t=0}= u(T/2) \quad\mbox{on}\;\;\R^3.\\
\end{array}\right.
\leqno(NS_{\nu,T}) $$ By Theorem \ref{t02}, there exists a
unique $v\in\mathcal C([0,T^*),H^s)$ solution of $(NS_{\nu,T})$,
satisfying
$$T^*<\infty\Longrightarrow \int_0^{T^*}\|v\|_{\dot H^s}=+\infty.$$
Suppose that $T^*\leq T_0-\frac{T}{2}$. By uniqueness we have $$\forall t\in[T/2,T^*),\;\;\;
u(t)=v(t-\frac{T}{2}).$$Taking the scalar product in $\dot H^s$,
and using lemma \ref{l1}, we obtain $\forall t\in[0,T^*)$
$$\begin{array}{lcl}
\partial_t\|v(t)\|_{\dot H^s}^2+2\nu\|\nabla v(t)\|_{\dot
H^s}^2&\leq& C\|\nabla v(t)\|_{L^2}\|v(t)\|_{\dot
H^s}^{1/2}\|\nabla v(t)\|_{\dot H^s}^{3/2}\\\\
&\leq& C\|\nabla v(t)\|_{L^2}^4\|v(t)\|_{\dot
H^s}^2+\nu\|\nabla v(t)\|_{\dot H^s}^{2}.\\
\end{array}$$
By Gronwall's lemma
$$\begin{array}{lcl}
\|v(t)\|_{\dot H^s}^2 &\leq& \|v(0)\|_{\dot
H^s}^2\exp\Big({C\displaystyle\int_0^t\|\nabla v(t)\|_{L^2}^4}\Big).\\\\
&\leq& \|u(T/2)\|_{\dot
H^s}^2\exp\Big({C\displaystyle\int_{T/2}^{t+\frac{T}{2}}\|\nabla u\|_{L^2}^4}\Big).\\\\
&\leq& \|u(T/2)\|_{\dot
H^s}^2\exp\Big({C\|u\|_{L^\infty_{T_0}(\dot H^{1/2})}^2\|u\|_{L^2_{T_0}(\dot H^{3/2})}^2}\Big).\\\\
\end{array}$$
Then $T^*>T_0-\frac{T}{2}$, we obtain $(\beta_1')$ consequently $(\beta_3')$.\\
Combines the first step and Theorems \ref{t1} we obtain
$(\beta_2')$. This completes the proof.
\endproof
\subsection{Proof Of Theorem \ref{apt2}} The proprieties $(\alpha_1')$ and
$(\alpha_2')$ are an easy consequences of Theorem \ref{apt1}.\\
{\bf Proof of $(\alpha_3')-(\alpha_4')$ :} If we prove the existence
of a time $T\geq 0$ such that $\|u(T)\|_{\dot H^{1/2}}\leq c\nu$, we
can apply Theorem \ref{t1} on $[T,+\infty)$, by the uniqueness we
obtain the desired results. Then, for simplification, we begin by
proving the following assertion.
\begin{equation}
\label{p1}
 \forall\; t\geq
0, \quad\|u\|_{L^\infty_t(L^2)}\leq \|u^0\|_{L^2}.
\end{equation}

Let $t^*:=\sup\{t\geq 0,\;\;\ \|u\|_{L^\infty_t(L^2)}\leq
\|u^0\|_{L^2}\}\in [0,+\infty]$. Suppose that $t^*<\infty$, by
continuity of $u$ show that $\|u(t^*)\|_{L^2}=\|u^0\|_{L^2}$.
Applying Theorem \ref{t1} to the following system
$$ \left\{\begin{array}{l}
  \displaystyle\partial_t
v-  \nu\Delta v+ v.\nabla v=-\nabla p,\quad
\mbox{on}\; \R_+\times\R^3, \\\\
div\; (v) = 0 \quad
\mbox{on}\; \R_+\times\R^3,  \\\\
v_{\mid t=0}= u(t^*) \quad\mbox{on}\;\R^3.\\
\end{array}\right.
\leqno(NS_{\nu}^{t^*}) $$ we obtain a time $T_1>0$ and a unique
solution $v\in\mathcal C([0,T_1],H^{1/2})\cap L^2_{T_1}(H^{3/2})$, satisfying
$$\forall t\in[0,T_1],\;\;\;
\|v(t)\|_{L^2}^2+2\nu\|\nabla v\|_{L^2_t(L^2)}^2\leq
\|v(0)\|_{L^2}^2=\|u^0\|_{L^2}^2.
$$
Using the uniqueness, we obtain $$u(t)=v(t-t^*),\;\;\forall
t\in[t^*,t^*+T_1],$$ then $$\|u(t)\|_{L^2}\leq
\|u^0\|_{L^2},\;\;\;\forall t\in[0,t^*+T_1].$$ Hence $t^*=+\infty$ and the assertion \eqref{p1} is proved.\\
Now, we have to prove
\begin{equation}
\label{p2}
\exists\; t\geq 0\quad\mbox{s.t.}\quad \|u(t)\|_{\dot H^{1/2}}\leq
c\nu.
\end{equation}
Let $$A:=\{t\geq 0,\;\;\ \|u(t)\|_{\dot H^{1/2}}> c\nu\}.$$ Using
H\"older inequality, we infer
$$\forall t\geq 0,\;\;\; (c\nu)^6{\bf
1}_A(t)\leq\|u(t)\|_{\dot H^{1/2}}^6\leq
\|u(t)\|_{L^2}^4\|u(t)\|_{\dot H^{3/2}}^2\leq
\|u^0\|_{L^2}^4\|u(t)\|_{\dot H^{3/2}}^2.
$$
We integrate on $\R^+$, we obtain
$$\lambda_1(A)\leq \frac{\|u^0\|_{L^2}^4\|u\|_{L^2(\dot
H^{3/2})}^2}{(c\nu)^6}:=t_0<+\infty,$$
where $\lambda_1$ is the Lebesgue measure on $\R$.\\
Then, for $\mu>0$, there exist $t_\mu\in(0,t_0+\mu)$, such that
$\|u(t_\mu)\|_{\dot H^{1/2}}\leq
c\nu$.\\Now, we consider the following
system
$$ \left\{\begin{array}{l}
  \displaystyle\partial_t
v-  \nu\Delta v+ v.\nabla v=-\nabla p,\quad
\mbox{on}\; \;\R_+\times\R^3, \\\\
div\; (v) = 0 \quad
\mbox{on}\; \;\R_+\times\R^3,  \\\\
v_{\mid t=0}= u(t_1) \quad\mbox{on}\;\;\R^3.\\
\end{array}\right.
\leqno(NS_{\nu}^{t_1}) $$ By Theorem \ref{t1} there exists a
unique $v_1\in\mathcal C_b(\R^+,H^{1/2})\cap L^2(\R^+,\dot
H^{3/2})$, solution of $(NS_{\nu}^{t_1})$ satisfying
$(\alpha_3)-(\alpha_4)$. The uniqueness imply $u(t)=v_1(t-t_1)$
for all $t\in[t_1,+\infty)$. This completes the proof.
\endproof

\end{document}